\documentclass[12pt]{article}
\usepackage{amsfonts, amssymb, amsmath}
\usepackage{enumerate}
\usepackage{esint}

\newtheorem{thm}{Theorem}[section]

\newtheorem{prop}[thm]{Proposition}
\newtheorem{defn}[thm]{Definition}
\newtheorem{rem}[thm]{Remark}

\newcommand{\f}{\frac}

\newcommand{\vc}{\infty}

\newcommand{\RR}{\mathbb{R}^n}

\newcommand{\p}{\partial}
\begin{document}

\title{Weighted norm inequalities for multilinear operators and applications to multilinear Fourier multipliers}

\author{The Anh Bui\thanks{The Anh Bui was supported by a Macquarie University scholarship} \and
Xuan Thinh Duong
\thanks{Xuan Thinh Duong was supported by a research grant from Macquarie University and Australian
Research Council (ARC) \newline
{\it {\rm 2010} Mathematics Subject Classification: 42B15, 42B20}.
}}

\date{}

\maketitle
\begin{abstract}
Let $T$ be a multilinear operator which is bounded on certain products of unweighted Lebesgue spaces of $\mathbb R^n$.
We assume that the associated kernel of $T$ satisfies some mild regularity condition  which is weaker than the usual
H\"older continuity   of those in the class of multilinear Calder\'on-Zygmund singular integral operators.
We then show the boundedness for $T$ and the boundedness of the commutator of $T$ with BMO functions
on products of weighted Lebesgue spaces of $\mathbb R^n$. As an application, we obtain the weighted norm inequalities
of  multilinear Fourier multipliers and of their commutators with BMO functions on the products of weighted Lebesgue spaces
when the number of derivatives of the symbols is the same as the best known result for the multilinear Fourier multipliers
to be bounded on the products of unweighted Lebesgue spaces.
\end{abstract}
\tableofcontents

\section{Introduction}
The theory of multilinear Calder\'on-Zygmund singular integral operators, originated from
the work of Coifman and Meyer, has had an important role in harmonic analysis. This direction of research has
been attracting a lot of attention in the last few decades, see for example \cite{CM1,CM2,CM3,GT, KS, LOPTG}
for the standard theory of multilinear Calder\'on-Zygmund singular integrals.
Recently, there are a number of studies concerning multilinear singular integrals which possess rough associated kernels
so that they do not belong to the standard  Calder\'on-Zygmund classes. See, for example
\cite{DGY, DGGLY, GLY, T, GS} and the references therein. \\

In this paper, we aim to study the boundedness of
multilinear singular integral operators on product of weighted Lebesgue spaces.
Assume that $T$ is a multilinear operator initially defined on the $m$-fold
product of Schwartz spaces and taking values into the space of
tempered distributions,
\begin{equation*}
T: \mathcal{S}(\mathbb{R}^n)\times\ldots\times
\mathcal{S}(\mathbb{R}^n) \rightarrow \mathcal{S}'(\mathbb{R}^n)
\end{equation*}
By the associated kernel $K(x, y_1, \ldots, y_m)$, we mean that $K$ is a function defined off the diagonal $x = y_1 =\ldots=y_m$ in
$(\mathbb{R}^n)^{m+1}$, satisfying
\begin{equation*}
T (f_1, \cdots, f_m)(x) = \int_{(\mathbb{R}^n)^m} K(x,y_1,\ldots ,
y_m)f_1(y_1)\ldots f_m(y_m)dy_1 \ldots dy_m
\end{equation*}
for   all functions $f_j\in \mathcal{S}(\mathbb{R}^n)$ and all $x \notin \cap^m_{j=1}$supp$f_j$, $j=1,\ldots, m$.

\medskip

In what follows, we denote $dy_1\ldots dy_m$ by $d\vec{y}$. For the rest of this paper, we assume that there exist $p_0\geq 1$ and a constant $C>0$ so that the following conditions holds:

\begin{enumerate}[\textbf{(H1)}]
\item \ \ \  For all $p_0\leq  q_1,  q_2,  \ldots,  q_m<\vc $  and $0 < q < \vc$ with
$$
\f{1}{q_1}+\ldots+\f{1}{q_m}=\f{1}{q},
$$
$T$ maps $L^{q_1}\times\ldots\times L^{q_m}$ into $L^{q,\vc}$.
\end{enumerate}
\begin{enumerate}[\textbf{(H2)}]
\item \ \ \ There exists $\delta> n/p_0$ so that for the conjugate exponent $p'_0$ of $p_0$, one has
\begin{equation}\label{cond2}
\begin{aligned}
\Big(\int_{S_{j_m}(Q)}\ldots \int_{S_{j_1}(Q)}&|(K(x,y_1, y_m)-K(\overline{x},y_1,\ldots, y_m))|^{p'_0}d\vec{y}\Big)^{1/{p'_0}}\\
&\leq C\f{|x-\overline{x}|^{m(\delta-n/{p_0})}}{|Q|^{m\delta/n}} 2^{-m\delta j_0}
\end{aligned}
\end{equation}
for all balls $Q$, all $x, \overline{x}\in \f{1}{2}Q$ and $(j_1,\ldots,y_m)\neq (0,\ldots,0)$, where $j_0=\max\{j_k: k=1,\ldots, m\}$ and $S_j(B)=2^{j}Q\backslash 2^{j-1}Q$ if $j\geq 1$, otherwise $S_j(Q)=Q$.
\end{enumerate}

Under the assumptions  {\bf{(H1)} } and \textbf{(H2)}, the multilinear operator $T$ may not fall under the scope of the theory of multilinear Calder\'on-Zygmund singular integral operators in \cite{GT}. An important example is the multilinear Fourier multiplier operators:
$$
T_m(f,g)(x)=\f{1}{(2\pi)^{2n}}\int_{\RR}\int_{\RR}e^{ix\cdot (\xi+\eta)}m(\xi,\eta)\hat{f}(\xi)\hat{g}(\eta)d\xi d\eta
$$
for all $f, g\in \mathcal{S}(\RR)$ when the function $m$ is not sufficiently smooth (see Section 3 for precise definition).
Indeed, if the number of derivatives imposed on the
function $m$ is not large enough, one may not expect the standard pointwise estimate for the kernel of $T_m$ in general.
The aim of this paper is to prove the weighted norm inequalities of such a multilinear operator $T$ and
weighted estimates of the commutator of $T$ with a BMO function.
We then consider the  multilinear Fourier multiplier operator $T = T_m$
in which the number of derivatives of the symbol $m$ is the same as that of \cite{T, GS} which guarantees the multilinear Fourier multiplier $T_m$
to be bounded on the products of unweighted Lebesgue spaces  \cite[Corollary 1.2]{T} and \cite[Theorem 1.1]{ GS}.
In this setting, we obtain  the weighted norm inequalities for $T_m$ and for the commutators of $T_m$ and a BMO function.\\

The layout of the paper is as follows. In Section 2, we recall some basic properties on weighted estimates for some maximal operators. The results on weighted estimates of multilinear operators and their commutators with BMO functions will be addressed in Section 3. As an application, we
 will consider in Section 4 the weighted norm inequalities for multilinear Fourier multiplier operators and their commutators with BMO functions.
\section{Sharp maximal function and weighted estimates}
\subsection{Sharp maximal operators}
We denote the Hardy-Littlewood maximal function with respect to balls on $\mathbb{R}^n$
by $M$. For $\delta>0$, let $M_{\delta}$ be the maximal function
\begin{equation*}
M_\delta
f(x)=M(|f|^\delta)^{1/\delta}=\Big(\sup\limits_{Q\ni x}\frac{1}{|Q|}\int_{Q}|f(y)|^\delta
dy\Big)^{1/\delta}.
\end{equation*}
Also, let $M^\sharp$ be the standard sharp maximal function of
Fefferman and Stein,
\begin{equation*}
M^\sharp
f(x)=\sup\limits_{Q\ni x} \ \inf_{c}\frac{1}{|Q|}\int_{Q}|f(y)-c|dy \ \approx \
\sup\limits_{Q\ni x}\frac{1}{|Q|}\int_{Q}|f(y)-f_Q|dy,
\end{equation*}
where $f_Q=\frac{1}{|Q|}\int_{Q}f(y)dy$ and $M^\sharp_\delta$
is defined by $M^\sharp_\delta(f)=(M^\sharp(|f|^\delta))^{1/\delta}$.\\

\medskip

We will denote the Muckenhoupt class by $A_\infty$. Let $\omega$ be a weight in
the Muckenhoupt class $A_\infty$ and let $0<p, \delta<\infty$. It is well known that (see, for example \cite{FS})  there exists $C>0$ (depending on the $A_\infty$
constant of $\omega$) such that
\begin{equation}\label{sharpweightinequality}
\int_{\mathbb{R}^n}(M_\delta f(x))^p\omega(x)dx\leq C
\int_{\mathbb{R}^n}(M^\sharp_{\delta} f(x))^p\omega(x)dx,
\end{equation}
for any function $f$ for which the left hand side is finite.\\

\subsection{Multiple weights}

For $m$ exponents $p_1, . . . , p_m$, we denote by $p$  the number given by
$1/p=1/p_1+\ldots+1/p_m$, and $\vec{P}$ for the vector $\vec{P} = (p_1, \ldots , p_m)$. The following definition of the class of multiple weights $A_{\vec{P}}$ is taken from \cite{LOPTG}.
\begin{defn}
Let $1\leq p_1, . . . , p_m<\infty$. Given $\vec{\omega} = (\omega_1, \ldots , \omega_m)$, set
\begin{equation*}
v_{\vec{\omega}}=\prod_{j=1}^{m}\omega_j^{p/p_j}
\end{equation*}
for all balls $Q$.
We say that $\vec{\omega}$ satisfies the $A_{\vec{P}}$ condition if
\begin{equation*}
\sup_{Q}\Big(\frac{1}{|Q|}\int_{Q}v_{\vec{\omega}}\Big)^{1/p}\prod^{m}_{j=1}\Big(\frac{1}{|Q|}\int_{Q}\omega_j^{1-p'_j}\Big)^{1/p'_j}<\infty.
\end{equation*}
When $p_j=1, \Big(\frac{1}{|Q|}\int_{Q}\omega_j^{1-p'_j}\Big)^{1/p'_j}$ is understood as $(\inf_{Q}\omega_j)^{-1}$.
\end{defn}

\begin{rem}\label{rem1}
{\rm Note that if $\vec{\omega}\in A_{\vec{P}}$ then $v_{\vec{\omega}}\in A_{mp}$ and there exists $\min\{p_1,\ldots,p_m\}>r>1$ such that $\vec{\omega}\in A_{\vec{P}/r}$, where $\vec{P}/r= (p_1/r, \ldots , p_m/r)$, see \cite{LOPTG}}.
\end{rem}
For $\vec{f}=(f_1,\ldots, f_m)$ and $p\geq 1$ we define the operator $\mathcal{M}_p$ by setting
$$
\mathcal{M}_p(\vec{f})(x)=\sup_{Q\ni x}\prod\limits_{j=1}^m \Big(\f{1}{|Q|}\int_Q|f_j(y_j)|^pdy_j\Big)^{1/p}.
$$
Note that the operator $\mathcal{M}_p$ when $p=1$ was introduced by \cite{LOPTG}. When $p=1$, we write $\mathcal{M}$ instead of $\mathcal{M}_1$. We have the following the weighted estimate for multilinear operator $\mathcal{M}_p(\vec{f})$.
\begin{prop}\label{weightine1}
Let $p_0\geq 1$ and  $p_j>p_0$ for all $j=1,\ldots,m$ and $\frac{1}{p}=\frac{1}{p_1}+\ldots+\frac{1}{p_m}$. Then
\begin{equation*}
\Big\|\mathcal{M}_{p_0}(\vec{f})\Big|\Big|_{L^p(v_{\vec{\omega}})}\leq C\prod_{j=1}^{m}||f_j||_{L^{p_j}(\omega_j)}
\end{equation*}
if and only if $\vec{\omega}\in A_{\vec{P}/p_0}$, where $\vec{P}/p_0=(p_1/p_0,\ldots,p_m/p_0)$.
\end{prop}
\emph{Proof:} We have $\mathcal{M}_{p_0}(\vec{f})=\mathcal{M}(\vec{f}^{p_0})^{1/p_0}$ where $\vec{f}^{p_0}=(f_1^{p_0}, \ldots, f_m^{p_0})$. Hence
$$
\Big\|\mathcal{M}_{p_0}(\vec{f})\Big\|_{L^p(v_{\vec{\omega}})}=\Big\|\mathcal{M}(\vec{f})\Big\|^{p_0}_{L^{p/p_0}(v_{\vec{\omega}})}.
$$
Using \cite[Theorem 3.7]{LOPTG}, we obtain
$$
\Big\|\mathcal{M}(\vec{f})\Big\|^{p_0}_{L^{p/p_0}(v_{\vec{\omega}})}\leq C\prod\limits_{j=1}^m \|f_j\|^{p_0}_{L^{p_j/p_0}(w_j)}=C\prod\limits_{j=1}^m \|f_j\|_{L^{p_j}(w_j)}.
$$
This completes our proof.
\section{Main results}
\subsection{Weighted estimates for multilinear operators}

The following result is an estimate on the Fefferman-Stein maximal function
acting on $T(\vec{f})$ in terms of Hardy-Littlewood maximal function.

\begin{thm}\label{mainthm1}
Let $T$ satisfy ($H_1$) and ($H_2$) and let $0<\delta<p_0/m$. Then for any $\vec{f}$ in the product space
$L^{q_1}(\mathbb{R}^n) \times L^{q_2}(\mathbb{R}^n) \times \cdots \times L^{q_m}(\mathbb{R}^n)$
with $p_0\leq q_j\leq \infty$ for $j = 1, 2, \cdots m$, one has
\begin{equation*}
M_\delta^\sharp(T(\vec{f}))(x)\leq
C\mathcal{M}_{p_0}(\vec{f})(x).
\end{equation*}
\end{thm}
\emph{Proof:} Fix a point $x$ and a ball $Q\ni x$. Due to the fact
that $\Big||\alpha|^r-|\beta|^r\Big|\leq |\alpha-\beta|^r$ for all $0<r<1$, we need only to prove that
\begin{equation*}
\Big(\frac{1}{|Q|}\int_{Q}|T(\vec{f})(z)-c_Q|^{\delta}dz\Big)^{1/\delta}\leq
C\mathcal{M}_{p_0}(\vec{f})(x),
\end{equation*}
where the constant $c_Q$ is to be fixed later and depends on $Q$. Using the standard argument, see for example \cite{GT, LOPTG}, for each $j$ we decompose $f_j=f_j^0+f_j^\infty$, where
$f_j^0=f_j\chi_{Q^*}, j=1,\ldots, m,$ and $Q^*=8Q$.
Then
\begin{equation*}
\begin{aligned}
\prod_{j=1}^{m}f_j(y_j)&=\sum_{\alpha_1,\ldots,\alpha_m\in\{0,\infty\}}f_1^{\alpha_1}(y_1)\ldots
f_m^{\alpha_m}(y_m)\\
&=\prod_{j=1}^{m}f_j^0(y_j)+ \sum_{(\alpha_1,\ldots,\alpha_m)\in \mathcal{I}_{\alpha}}f_1^{\alpha_1}(y_1)\ldots
f_m^{\alpha_m}(y_m)
\end{aligned}
\end{equation*}
where $\mathcal{I}_{\alpha}=\{(\alpha_1,\ldots,\alpha_m): \ \text{there is at least one $\alpha_j\neq 0$}\}$. So, we can  write
\begin{equation*}
T(\vec{f})(z)=T(\vec{f}^0)(z)+\sum_{(\alpha_1,\ldots,\alpha_m)\in \mathcal{I}_{\alpha}}T(f_1^{\alpha_1}\ldots
f_m^{\alpha_m})(z).
\end{equation*}
Due to (H1), $T$ maps $L^{p_0}\times \ldots \times L^{p_0}$ into $L^{{p_0}/m,\infty}$. This together with Kolmogorov inequality tells us that
\begin{equation*}
\begin{aligned}
\Big(\frac{1}{|Q|}\int_{Q}|T(\vec{f}^0)(z)|^\delta
dz\Big)^{1/\delta}&\leq
C||T(\vec{f}^0)(z)||_{L^{p_0/m,\infty}(Q,\frac{dx}{|Q|})}\\
&\leq C\prod_{j=1}^{m}\Big(\frac{1}{|Q^*|}\int_{Q^*}|f_j(z)|^{p_0}dz\Big)^{p_0}\\
&\leq C \mathcal{M}_{p_0}(\vec{f})(x).
\end{aligned}
\end{equation*}
To estimate the remaining terms, we choose $c_Q=\sum_{(\alpha_1,\ldots,\alpha_m)\in \mathcal{I}_{\alpha}}(f_1^{\alpha_1}\ldots f_m^{\alpha_m})(x)$. We will verify that
\begin{equation}\label{equation1}
\sum_{(\alpha_1,\ldots,\alpha_m)\in \mathcal{I}_{\alpha}}|T(f_1^{\alpha_1}\ldots
f_m^{\alpha_m})(z)-T(f_1^{\alpha_1}\ldots f_m^{\alpha_m})(x)|\leq C
\mathcal{M}_{p_0}(\vec{f})(x).
\end{equation}
For $(\alpha_1,\ldots,\alpha_m)\in \mathcal{I}_{\alpha}$, we assume that $\alpha_1=\ldots=\alpha_l =\vc$ and $\alpha_{l+1}=\ldots=\alpha_m =0$, $l\geq 1$. For such a $(\alpha_1,\ldots,\alpha_m)$, we can write
\begin{equation*}
\begin{aligned}
|T(f_1^{\alpha_1}&\ldots f_m^{\alpha_m})(z)-T(f_1^{\alpha_1}\ldots
f_m^{\alpha_m})(x)|\\
&\leq
\int_{(\mathbb{R}^n)^m}|K(z,\vec{y})-K(x,\vec{y})\prod\limits_{j=1}^m|f_j(y_j)|d\vec{y}\\
&\leq  \int_{(\mathbb{R}^n\backslash Q^*)^l\times (Q^*)^{m-l}}|K(z,\vec{y})-K(x,\vec{y})|\prod\limits_{j=1}^l|f^\vc_j(y_j)|\prod\limits_{j=l+1}^m|f^0_j(y_j)|d\vec{y}\\
&\leq  \sum_{j_1,...,j_l\geq 1}\int_{(Q^*)^{m-l}}\int_{S_{j_l}(Q^*)}\ldots\int_{S_{j_1}(Q^*)}|K(z,\vec{y})-K(x,\vec{y})|\prod\limits_{j=1}^l|f^\vc_j(y_j)|\prod\limits_{j=l+1}^m|f^0_j(y_j)|d\vec{y}\\
\end{aligned}
\end{equation*}
Using H\"older inequality and (H2), we have
 \begin{equation*}
\begin{aligned}
\sum_{j_1,...,j_l\geq 1}\int_{(Q^*)^{m-l}}&\int_{S_{j_l}(Q^*)}\ldots\int_{S_{j_1}(Q^*)}|K(z,\vec{y})-K(x,\vec{y})|\prod\limits_{j=1}^l|f^\vc_j(y_j)|\prod\limits_{j=l+1}^m|f^0_j(y_j)|d\vec{y}\\
&\leq C\sum_{j_1,...,j_l\geq 1}\Big(\int_{(Q^*)^{m-l}}\int_{S_{j_l}(Q^*)}\ldots\int_{S_{j_1}(Q^*)}|K(z,\vec{y})-K(x,\vec{y})|^{p_0'}d\vec{y}\Big)^{1/p_0'}\\
&~~~\times \prod\limits_{j=1}^l\Big(\int_{2^{j_k}Q^*}|f_j(y_{j})|^{p_0}dy_j\Big)^{1/p_0}\prod\limits_{j=l+1}^m\Big(\int_{Q_*}|f_j(y_j)|^{p_0}dy_j\Big)^{1/p_0}\\
&\leq C\sum_{j_1,...,j_l\geq 1}\f{|x-\overline{x}|^{m(\delta-n/{p_0})}}{|Q^*|^{m\delta/n}} 2^{-m\delta j_0}\\
&~~~~\times\prod\limits_{j=1}^l\Big(\int_{2^{j_k}Q^*}|f_j(y_{j})|^{p_0}dy_j\Big)^{1/p_0}\prod\limits_{j=l+1}^m\Big(\int_{Q_*}|f_j(y_j)|^{p_0}dy_j\Big)^{1/p_0}\\
&\leq C\sum_{j_0\geq 1}\f{|x-z|^{m(\delta-n/{p_0})}}{|Q^*|^{m\delta/n}} m2^{-m\delta j_0} 2^{j_0mn/p_0}|Q^*|^{m/p_0}\\
&~~~~\times\prod\limits_{j=1}^m\Big(\f{1}{|2^{j_0}Q^*|}\int_{2^{j_0}Q^*}|f_j(y_{j})|^{p_0}dy_j\Big)^{1/p_0}\\
&\leq C\sum_{j_0\geq 1}\f{|x-z|^{m(\delta-n/{p_0})}}{|Q^*|^{m(\delta/n-1/p_0)}} m2^{-mj_0(\delta-n/p_0)}\mathcal{M}_{p_0}\vec{f}(x)\\
&\leq C\mathcal{M}_{p_0}\vec{f}(x)
\end{aligned}
\end{equation*}
as long as $\delta>n/p_0$ and $x, z\in Q$, where $j_0=\max\{j_1,....,j_l\}$.
This completes our proof.

\medskip

 The following theorem is our main result of weighted estimates for multilinear operators with rough kernels.

\begin{thm}\label{thmofT}
Let $T$ satisfy ($H_1$) and ($H_2$). For any $p_0< p_1, \ldots, p_m < \infty$ and $p$ so that $1/p_1+\ldots+1/p_m=1/p$ and $\vec{\omega}\in A_{\vec{P}/p_0}$, we have
\begin{equation*}
\|T(\vec{f})\|_{L^p(v_{\vec{\omega}})}\leq C\prod\limits_{j=1}^{m}\|f_j\|_{L^{p_j}(w_j)}.
\end{equation*}
\end{thm}

\emph{Proof:} The proof is just the combination of the results of Theorem \ref{mainthm1}, Proposition \ref{weightine1} and the weighted norm inequality (\ref{sharpweightinequality}).

\subsection{Weighted estimates for commutators of multilinear operators with BMO functions}
We next obtain an estimate  on the sharp maximal function of the commutator of a BMO vector function
and the multilinear operator. Given a locally integrable vector function
$\vec{b} = (b_1, . . . , b_m)$, we define the $m$-linear
commutator of $\vec{b}$ and the $m$-linear operator $T$ by
\begin{equation*}
T_{\vec{b}}(\vec{f})=\sum_{j=1}^{m}T_{\vec{b}}
^j(\vec{f})
\end{equation*}
where
\begin{equation*}
T_{\vec{b}}^j(\vec{f})=b_jT(\vec{f})-T(f_1,\ldots,
b_jf_j,\ldots, f_m).
\end{equation*}
We use the notation $||\vec{b}||_{BMO}=\max_j
||b_j||_{BMO}$.
\begin{thm}\label{thm2}
Assume that $T$ satisfies (H1) and (H2). Let $T_{\vec{b}}$ be a multilinear commutator with
$\vec{b} \in BMO^m$ and let  $0 < \delta < \epsilon <p_0/m$. Then for any $q_0>p_0$ there exists a constant $C >0$, depending
on $\delta$ and $\epsilon$, such that \begin{equation*}
M^\sharp_\delta(T_{\vec{b}}(\vec{f}))(x)\leq
C||\vec{b}||_{BMO}\mathcal{M}_{q_0}(\vec{f})(x)
\end{equation*}
for all bounded measurable
vector functions $\vec{f}=(f_1,\ldots, f_m)$ with compact supports.
\end{thm}
\emph{Proof:} By linearity it is sufficient to consider the particular case when $\vec{b}=b\in BMO$. Fix $b\in BMO$ and consider the operator
\begin{equation*}
T_b(\vec{f})(x)=b(x)T(\vec{f})-T(bf_1,\ldots,
f_m).
\end{equation*}
Fix $x\in \mathbb{R}^n$. For any ball $Q$ with center at $x$, set  $Q^*=8Q$. Then we have
\begin{equation*}
T_b(\vec{f})(x)=(b(x)-b_{Q^*})T(\vec{f})(x)-T((b-b_{Q^*})f_1,\ldots,
f_m)(x).
\end{equation*}
Since $0<\delta<1$,
\begin{equation*}
\begin{aligned}
\Big(\frac{1}{|Q|}\int_Q\Big||T_b(\vec{f})(z)|^\delta-|c|^\delta\Big|dz\Big)^{1/\delta}
&\leq
\Big(\frac{1}{|Q|}\int_Q|T_b(\vec{f})(z)-c|^\delta dz\Big)^{1/\delta}\\
&\leq  \Big(\frac{C}{|Q|}\int_Q|(b(z)-b_{Q^*})T(\vec{f})(z)|^\delta dz\Big)^{1/\delta}\\
& ~~+ \Big(\frac{C}{|Q|}\int_Q|T((b-b_{Q^*})f_1,\ldots,
f_m)(z)-c|^\delta dz\Big)^{1/\delta}\\
& = I+II
\end{aligned}
\end{equation*}
Pick $p>1$ so that $\delta p < \epsilon <p_0/m$ and $\delta p'>1$. By John-Nirenberg inequality and H\"older inequality, one has
\begin{equation*}
\begin{aligned}
I&\leq C\Big(\frac{1}{|Q|}\int_Q|(b(z)-b_{Q^*})|^{p'\delta} dz\Big)^{1/p'\delta}\Big(\frac{1}{|Q|}\int_Q|T(\vec{f})(z)|^{p\delta} dz\Big)^{1/p\delta}\\
&\leq C||b||_{BMO}\Big(\frac{1}{|Q|}\int_Q|T(\vec{f})(z)|^{p\delta} dz\Big)^{1/p\delta}\\
&\leq C||b||_{BMO}M_{\epsilon}(T(\vec{f}))(x) .
\end{aligned}
\end{equation*}
Using the similar decomposition to that in the proof of Theorem \ref{mainthm1}, we can write
\begin{equation*}
\begin{aligned}
\prod_{j=1}^{m}f_j(y_j)&=\sum_{\alpha_1,\ldots,\alpha_m\in\{0,\infty\}}f_1^{\alpha_1}(y_1)\ldots
f_m^{\alpha_m}(y_m)\\
&=\prod_{j=1}^{m}f_j^0(y_j)+ \sum_{(\alpha_1,\ldots,\alpha_m)\in \mathcal{I}_{\alpha}}f_1^{\alpha_1}(y_1)\ldots
f_m^{\alpha_m}(y_m).
\end{aligned}
\end{equation*}
Let $c=\sum_{(\alpha_1,\ldots,\alpha_m)\in \mathcal{I}_{\alpha}}T((b-b_{Q^*})f_1^{\alpha_1}\ldots f_m^{\alpha_m})(x)$.
We have
\begin{equation*}
\begin{aligned}
II&\leq C\Big(\Big(\frac{1}{|Q|}\int_Q|T((b-b_{Q^*})f_1^0,\ldots,
f_m^0)(z)|^\delta dz\Big)^{1/\delta}\\
&~~~+\sum_{(\alpha_1,\ldots,\alpha_m)\in \mathcal{I}_{\alpha}}\Big(\frac{1}{|Q|}\int_Q|T((b-b_{Q^*})f_1^{\alpha_1},\ldots,
f_m^{\alpha_m})(z)-T((b-b_{Q^*})f_1^{\alpha_1},\ldots,
f_m^{\alpha_m})(x)|^\delta dz\Big)^{1/\delta} \\
&=II_1+\sum_{(\alpha_1,\ldots,\alpha_m)\in \mathcal{I}_{\alpha}} II_{\alpha_1\ldots\alpha_m}.
\end{aligned}
\end{equation*}
We estimate the term $II_1$  by using
Kolmogorov inequality and H\"older inequality,
\begin{equation*}
\begin{aligned}
II_1& \leq C\|T((b-b_{Q^*})f_1^0,\ldots,
f_m^0)\|_{L^{p_0/m,\infty}(Q,\frac{dx}{|Q|})}\\
&\leq C\Big(\frac{1}{|Q|}\int_Q|(b-b_{Q^*})f_1^0(z)|dz\Big)^{1/p_0}\prod_{j=2}^\infty \Big(\frac{1}{|Q|}\int_Q|f_j^0(z)|^{p_0}dz\Big)^{1/p_0}\\
&\leq C||b||_{BMO}\mathcal{M}_{q_0}(\vec{f})(x).
\end{aligned}
\end{equation*}
Concerning the second term $\sum\limits_{(\alpha_1,\ldots,\alpha_m)\in \mathcal{I}_\alpha}
II_{\alpha_1\ldots\alpha_m}$, by using a similar argument to that in the proof of Theorem \ref{mainthm1}, we obtain that
\begin{equation*}
\begin{aligned}
\sum\limits_{(\alpha_1,\ldots,\alpha_m)\in \mathcal{I}_\alpha}
II_{\alpha_1\ldots\alpha_m}&\leq C\sum_{j_0\geq 1}\f{|x-z|^{m(\delta-n/{p_0})}}{|Q^*|^{m\delta/n}} m2^{-m\delta j_0} 2^{j_0mn/p_0}|Q^*|^{m/p_0}\\
&~~~~\times\Big(\f{1}{|2^{j_0}Q^*|}\int_{2^{j_0}Q^*}|(b(y_1)-b_{Q^*})f_1(y_{1})|^{p_0}dy_1\Big)^{1/p_0}\\
&~~~~\times\prod\limits_{j=2}^m\Big(\f{1}{|2^{j_0}Q^*|}\int_{2^{j_0}Q^*}|f_j(y_{j})|^{p_0}dy_j\Big)^{1/p_0}\\
&\leq C\|b\|_{BMO}\mathcal{M}_{q_0}\vec{f}(x)
\end{aligned}
\end{equation*}
provided $\delta>n/p_0$ and $x, z\in Q$. This completes our proof.

\medskip

As a consequence of Theorem \ref{thm2}, we have the following result.
\begin{thm}\label{thmofbT}
Let $T$ satisfy ($H_1$) and ($H_2$), and let $\vec{b} \in BMO^m$. For any $p_0< p_1, \ldots, p_m < \infty$ and $p$ such that $1/p_1+\ldots+1/p_m=1/p$ and $\vec{\omega}\in A_{\vec{P}/p_0}$, we have
\begin{equation*}
\|T_{\vec{b}}(\vec{f})\|_{L^p(v_{\vec{\omega}})}\leq C\|\vec{b}\|_{BMO}\prod\limits_{j=1}^{m}\|f_j\|_{L^{p_j}(w_j)}.
\end{equation*}
\end{thm}
\emph{Proof:} Since $\vec{\omega}\in A_{\vec{P}/p_0}$, there exists $r>1$ so that $\vec{\omega}\in A_{\vec{P}/p_0r}$. Taking $q_0=rp_0>p_0$, by Theorem \ref{thm2}, we have
$$
\|T_{\vec{b}}(\vec{f})\|_{L^p(v_{\vec{\omega}})}\leq C\|\vec{b}\|_{BMO}\Big(\|M_{\epsilon}(T(\vec{f}))(x)\|_{L^p(v_{\vec{\omega}})}+\|\mathcal{M}_{q_0}(\vec{f})\|_{L^p(v_{\vec{\omega}})}\Big).
$$
Since $\vec{\omega}\in A_{\vec{P}/p_0}$ and $\vec{\omega}\in A_{\vec{P}/q_0}$, using Proposition \ref{weightine1}, we have
$$
\|\mathcal{M}_{q_0}(\vec{f})\|_{L^p(v_{\vec{\omega}})}\leq C\prod\limits_{j=1}^{m}\|f_j\|_{L^{p_j}(w_j)} \ .
$$
Moreover, Remark \ref{rem1} tells us that $v_{\vec{\omega}}\in A_{pm/p_0}$. This together with Theorem \ref{thmofT} gives
$$
\|M_\epsilon(T(\vec{f}))\|_{L^p(v_{\vec{\omega}})}\leq C\|T(\vec{f})\|_{L^p(v_{\vec{\omega}})} \leq C\prod\limits_{j=1}^{m}\|f_j\|_{L^{p_j}(w_j)}.
$$
This completes our proof.

\section{Application to multilinear Fourier multipliers}

In this section, we apply the results in Section 3 to investigate the weighted estimates for multilinear Fourier multiplier operators. Before coming to the details, we consider the linear case first. Let $m\in L^{\vc}(\mathbb{R}^n)$. The Fourier multiplier operator $T_m$ is defined by
$$
T_mf(x)=\f{1}{(2\pi)^n}\int_{\mathbb{R}^n}e^{ix\cdot \xi}m(\xi)\hat{f}(x)d\xi
$$
for all Schwart functions $f\in S(\mathbb{R}^n)$, where $\hat{f}$ is the Fourier transform of $f$.

It is well-known that if $m$ satisfies the following condition
$$
|\p_\xi^\alpha m(\xi)|\leq C_\alpha|\xi|^{-\alpha} \ \ \text{for all $\alpha\leq [n/2]+1$}
$$
then $T_m$ is bounded on $L^p$ for all $1<p<\vc$, see for example \cite[Corollary 8.11]{Du}.

We now consider the multilinear case. For the sake of simplicity, we only consider the bilinear case. Let $m\in C^{s}(\mathbb{R}^{2n}\backslash\{0\})$, for some integer $s$, satisfying the following condition:
\begin{equation}\label{Homandercondition-multipliers}
|\p^\alpha_{\xi}\p^\beta_{\eta}m(\xi,\eta)|\leq C_{\alpha,\beta}(|\xi|+|\eta|)^{-(|\alpha|+|\beta|)}
\end{equation}
for all $|\alpha|+|\beta|\leq s$ and $(\xi,\eta)\in \mathbb{R}^{2n}\backslash\{0\}$. The bilinear Fourier multiplier operator $T_m$ is defined by
$$
T_m(f,g)(x)=\f{1}{(2\pi)^{2n}}\int_{\RR}\int_{\RR}e^{ix\cdot (\xi+\eta)}m(\xi,\eta)\hat{f}(\xi)\hat{g}(\eta)d\xi d\eta
$$
for all $f, g\in \mathcal{S}(\RR)$. \\

Concerning the boundedness of $T_m$, it was proved in \cite{CM2} that if (\ref{Homandercondition-multipliers}) holds for $s>4n$ then $T_m$ maps from $L^{p_1}\times L^{p_2}$ into $L^p$ for all $1<p_1, p_2, p<\vc$ so that $1/p_1+1/p_2=1/p$. Then, in \cite{GT}, the authors proved that $T_m$ maps boundedly from $L^{p_1}\times L^{p_2}$ into $L^p$ for all $1<p_1, p_2<\vc$ so that $1/p_1+1/p_2=1/p$ provided that (\ref{Homandercondition-multipliers}) holds for $s\geq 2n+1$. However, in the sense of the linear case, the number of derivatives $s\geq 2n+1$  is not optimum and it is natural to expect that we only need $s\geq n+1$. The first positive answer is due to Tomita \cite{T} who proved that if (\ref{Homandercondition-multipliers}) holds for $s\geq n+1$, then $T_m$ maps from $L^{p_1}\times L^{p_2}$ into $L^p$ for all $2\leq p_1, p_2, p<\vc$ such that $1/p_1+1/p_2=1/p$ and then by using the multilinear interpolation and duality arguments, he obtained that $T_m$ maps from $L^{p_1}\times L^{p_2}$ into $L^p$ for all $1<p_1, p_2, p<\vc$ such that $1/p_1+1/p_2=1/p$. This result was then improved by \cite{GS} for $p\leq 1$ by using the $L^r$-based Sobolev space, $1<r\leq 2$. We would like to point out a particular case of the result in \cite[Theorem 1.1]{GS} in the following theorem.

\begin{thm}\label{thmofGS}
Assume that $(\ref{Homandercondition-multipliers})$ holds for some $n+1\leq s\leq 2n$. Then for any $p_1, p_2$ and $p$ such that $\f{2n}{s}<p_1,p_2<\vc$ and $1/p_1+1/p_2=1/p$, the operator $T_m$ maps from $L^{p_1}\times L^{p_2}$ into $L^p$.
\end{thm}

We remark that the number  $\f{2n}{s}$ in Theorem \ref{thmofGS} is contained implicitly in the proof of \cite[Theorem 1.1]{GS}.

It is natural to raise the question of weighted estimates for multilinear operators $T_m$ and their commutators with BMO functions. The aim of this section is to give a positive answer for this problem by using the results in Section 3. Our main results are formulated by the following theorem.

\begin{thm}\label{thmTm}
Assume that $(\ref{Homandercondition-multipliers})$ holds for some $n+1\leq s\leq 2n$. Then for any $p_1, p_2, p$ such that $r_0:=\f{2n}{s}<p_1,p_2<\vc$,
  $1/p_1+1/p_2=1/p$, and $\vec{\omega}=(w_1,w_2)\in A_{\vec{P}/r_0}$ with $\vec{P}=(p_1,p_2)$,  we have
\begin{enumerate}[(a)]
\item $
\|T_m(f_1, f_2)\|_{L^p(v_{\vec{\omega}})}\leq C\|f_1\|_{L^{p_1}(w_1)}\|f_2\|_{L^{p_2}(w_2)};
$

\item $\|(T_m)_{\vec{b}}(f_1, f_2)\|_{L^p(v_{\vec{\omega}})}\leq C\|\vec{b}\|_{BMO}\|f_1\|_{L^{p_1}(w_1)}\|f_2\|_{L^{p_2}(w_2)},$
for all $\vec{b}\in (BMO(\mathbb{R}^n))^2$.
\end{enumerate}
\end{thm}

It is easy to see that the associated kernel $K(x, y_1, y_2)$ to $T_m$ is given by
\begin{equation}\label{associatedkenerl}
K(x,y_1,y_2)= m\check{ } (x-y_1,x-y_2)
\end{equation}
where $m\check{ }$ is the inverse Fourier transform of $m$. We now show that the associated kernel $K$ satisfies (H2).
\begin{prop}\label{kernelcondtionH2}
For any $2\geq p> 2n/s$, we have,
\begin{equation}\label{eq1-kernelH2}
\Big(\int_{S_j(Q)}\int_{S_k(Q)}|K(x,y_1,y_2)-K(\overline{x},y_1,y_2)|^{p'}dy_1dy_2\Big)^{1/p'}\leq C\f{|x-\overline{x}|^{s-2n/{p}}}{|Q|^{s/n}} 2^{-s \max\{j,k\}}
\end{equation}
for all balls $Q$, all $x, \overline{x}\in \f{1}{2}Q$ and $(j,k)\neq (0,0)$.
\end{prop}
\emph{Proof:}
Due to (\ref{associatedkenerl}), we need only to show that
\begin{equation}\label{eq1-proofkernelH2}
\Big(\int_{S_j(Q)}\int_{S_k(Q)}|m\check{ } (x-y_1,x-y_2)-m\check{ } (\overline{x}-y_1,\overline{x}-y_2)|^{p'}dy_1dy_2\Big)^{1/p'}\leq C\f{|x-\overline{x}|^{s-2n/{p}}}{|Q|^{s/n}} 2^{-s \max\{j,k\}}
\end{equation}
for all balls $Q$, all $x, \overline{x}\in \f{1}{2}Q$ and $(j,k)\neq (0,0)$.

Using the change of variables, this is equivalent to that
\begin{equation}\label{eq1-proofkernelH2}
\Big(\int_{S_j(Q_{\overline{x}})}\int_{S_k(Q_{\overline{x}})}|m\check{ } (y+h,z+h)-m\check{ }(y,z)|^{p'}dydz\Big)^{1/p'}\leq C\f{|h|^{s-2n/{p}}}{|Q|^{s/n}} 2^{-s \max\{j,k\}}
\end{equation}
for $(j,k)\neq (0,0)$, where $h=x-\overline{x}$ and $Q_{\overline{x}}=Q-\overline{x}$.

Let $\Psi\in \mathcal{S}(\mathbb{R}^{2n})$ satisfying supp $\Psi\in \{(\xi,\eta): 1/2\leq |\xi|+|\eta|\leq 2\}$ and
$$
\sum_{j\in \mathbb{Z}} \Psi(2^{-j}\xi,2^{-j}\eta)= 1 \ \ \text{for all $(\xi,\eta)\in \mathbb{R}^{2n}\backslash \{0\}$}.
$$
Therefore, we can write
\begin{equation}\label{eq1-decompose}
m(\xi,\eta)= \sum_{j\in \mathbb{Z}} \Psi(2^{-j}\xi,2^{-j}\eta)m(\xi,\eta):=\sum_{j\in \mathbb{Z}}m_j(\xi,\eta)
\end{equation}
and hence supp $m_j\in \{(\xi,\eta): 2^{j-1}\leq |\xi|+|\eta|\leq 2^{j+1}\}$.

Without of the loss of generality, we assume that $k\geq j$ and hence $k\geq 1$. With the decomposition as in (\ref{eq1-decompose}), we set
$$
A_l =\Big(\int_{S_j(Q_{\overline{x}})}\int_{S_k(Q_{\overline{x}})}|m_l\check{ } (y+h,z+h)-m\check{ }_l (y,z)|^{p'}dydz\Big)^{1/p'}
$$
It is easy to see that  $2^{k-2}R\leq |y+h|\leq 2^{k+1}R$ and $|z+h|\leq 2^{j+1}R$ for all $y\in S_k(Q_{\overline{x}})$ and $z\in S_j(Q_{\overline{x}})$, where $R=l(Q)/2$. Therefore,
$$
A_l\leq C\Big(\int_{|z|<2^{j+1}R}\int_{2^{k-2}R\leq |y|\leq 2^{k+1}R}|m\check{ }_l (y,z)|^{p'}dydz\Big)^{1/p'}.
$$
Since $|x|\approx 2^kR$, by Hausdoff-Young inequality, we have, for $|\alpha|=s$,
\begin{equation*}
\begin{aligned}
A_l&\leq C(2^kR)^{-|\alpha|}\Big(\int_{|z|<2^{j+1}R}\int_{2^{k-2}R\leq |y|\leq 2^{k+1}R}|y|^{\alpha}|m\check{ }_l (y,z)|^{p'}dydz\Big)^{1/p'}\\
&\leq C2^{-ks}R^{-s}\sum_{|\alpha|=s}\Big(\int_{\RR}\int_{\RR}|\p_\xi^{\alpha}m_l (\xi,\eta)|^{p}d\xi d\eta\Big)^{1/p}
\end{aligned}
\end{equation*}
Using the fact that  supp $m_l\in \{(\xi,\eta): 2^{l-1}\leq |\xi|+|\eta|\leq 2^{l+1}\}$ and condition (\ref{Homandercondition-multipliers}), we get that
\begin{equation*}
A_l\leq C2^{-ks}R^{-s}2^{-l(s-2n/p)}.
\end{equation*}
This implies that
\begin{equation}\label{eq1-proofKerneles}
\sum_{l:2^{-l}\leq |h|}A_l \leq C2^{-ks}R^{-s}|h|^{s-2n/p}
\end{equation}
provided $2n/p <s$.

We can also write
$$
m\check{ }_l (y+h,z+h)-m\check{ }_l (y,z)= \varphi_l\check{ } (y,z)
$$
where $\varphi_l(\xi, \eta)=m_l(\xi, \eta)(e^{i(h\cdot \xi+h\cdot \eta)}-1)$.

Using Hausdoff-Young inequality again, we obtain that, for $|\alpha|=s$,
\begin{equation}\label{eq2-proofKerneles}
\begin{aligned}
A_l&\leq C(2^kR)^{-|\alpha|}\Big(\int_{S_j(Q_{\overline{x}})}\int_{S_k(Q_{\overline{x}})}|y|^{\alpha}|\varphi_l\check{ } (y,z)|^{p'}dydz\Big)^{1/p'}\\
&\leq C2^{-ks}R^{-s}\sum_{|\alpha|=s}\Big(\int_{\RR}\int_{\RR}|\p_\xi^{\alpha}\varphi_l (\xi,\eta)|^{p}d\xi d\eta\Big)^{1/p}\\
&\leq C2^{-ks}R^{-s}\sum_{|\alpha|=s}\Big(\int_{\RR}\int_{\RR}|\p_\xi^{\alpha}[m_l(\xi, \eta)(e^{i(h\cdot \xi+h\cdot\eta)}-1)]|^{p}d\xi d\eta\Big)^{1/p}.
\end{aligned}
\end{equation}
Moreover, we have
$$
\sum_{|\alpha|=s}|\p_\xi^{\alpha}[m_l(\xi, \eta)(e^{i(h\cdot \xi+h\cdot \eta)}-1)]|= \sum_{|\beta|+|\gamma|=s}|\p_\xi^\beta m_l(\xi, \eta)|\times |\p_\xi^\gamma(e^{i(h\cdot \xi+h\cdot \eta)}-1)|.
$$
Note that $|\p_\xi^\gamma(e^{i(h\cdot \xi+h\cdot \eta)}-1)|\leq |h||\xi|$ if $\gamma=0$, otherwise
$$|\p_\xi^\gamma(e^{i(h\cdot \xi+h\cdot \eta)}-1)|\leq |h|^{\gamma}\leq |h|2^{-l(|\gamma|-1)}$$
provided $2^{l}|h|<1$.

Therefore, for all $l$ with $2^l|h|<1$, we have
$$
\sum_{|\alpha|=s}|\p_\xi^{\alpha}[m_l(\xi, \eta)(e^{(i(h_1\cdot \xi+h_2\cdot \eta)}-1)]|\leq C|h_1|2^{-l(|\alpha|-1)}.
$$
This together with (\ref{eq2-proofKerneles}) gives
$$
A_l\leq C2^{-ks}R^{-s} |h|2^{-l(s-2n/p-1)}
$$
whenever $2^l|h_1|<1$.
Hence
\begin{equation}\label{eq3-proofKerneles}
\sum_{l:2^l|h|<1}A_l\leq C2^{-ks}R^{-s} |h|^{s-2n/p}
\end{equation}
as long as $2n/p<s$.

Combining (\ref{eq1-proofKerneles}) and (\ref{eq3-proofKerneles}), we complete the proof.\\

\emph{Proof of Theorem \ref{thmTm}:} Since $\vec{\omega}\in A_{\vec{P}/r_0}$, there exists $\min\{p_1/r_0, p_2/r_0\}>\alpha>1$ such that $\vec{\omega}\in A_{\vec{P}/\alpha r_0}$. Taking $p_0=\alpha r_0$, we have $p_1, p_2> p_0>r_0$. It follows from Theorem \ref{thmofGS} and Proposition \ref{kernelcondtionH2}  that $T_m$ satisfies (H1) and (H2) for $p_0$. Hence Theorem \ref{thmTm} is just a direct consequence of Theorems \ref{thmofT} and \ref{thmofbT}.

\medskip

\noindent Department of Mathematics, Macquarie University, NSW 2109, Australia and \\
Department of Mathematics, University of Pedagogy, Ho chi Minh city, Vietnam \\
Email: the.bui@mq.edu.au

\medskip

\noindent Department of Mathematics, Macquarie University, NSW 2109, Australia \\
Email: xuan.duong@mq.edu.au

\end{document}